\documentclass[a4paper,12pt]{amsart}

\usepackage{enumerate}
\usepackage{esint}
\usepackage[latin1]{inputenc}
\usepackage{amsmath}
\usepackage{amsthm}
\usepackage{latexsym}
\usepackage{amssymb}
\usepackage{amsmath}
\usepackage{comment}
\usepackage{calc}
\usepackage{multicol}
\usepackage{slashed}
\usepackage{graphicx}
\usepackage{bbm}
\newtheorem{thm}{Theorem}[section]
\newtheorem{prop}{Proposition}[section]

\newtheorem{cor}{Corollary}[section]

\newcommand{\R}{\mathbb{R}}
\numberwithin{equation}{section}

\newcommand{\Z}{\mathbb{Z}}

\newcommand{\vertiii}[1]{{\left\vert\kern-0.25ex\left\vert\kern-0.25ex\left\vert #1
\right\vert\kern-0.25ex\right\vert\kern-0.25ex\right\vert}}

\usepackage{xcolor}
\makeatletter
\newcommand{\leqnomode}{\tagsleft@true}
\newcommand{\reqnomode}{\tagsleft@false}
\makeatother

\begin{document}
\title[A note on the fractional Hardy inequality]{A note on the fractional Hardy inequality}
\author[M. Aldovardi and J. Bellazzini
]{Matteo Aldovardi and Jacopo Bellazzini }

\address{Matteo Aldovardi, Dipartimento di Matematica, Universit\`a Degli Studi di Pisa, Largo Bruno Pontecorvo, 5, 56127, Pisa, Italy}
\address{Jacopo Bellazzini, Dipartimento di Matematica, Universit\`a Degli Studi di Pisa, Largo Bruno Pontecorvo, 5, 56127, Pisa, Italy}


\thanks{J. B. is partially supported by project PRIN 2020XB3EFL by the Italian Ministry of Universities and Research and by  the University of Pisa, Project PRA 2022 11.}

\maketitle
\begin{abstract}
We give a direct proof of fractional Hardy inequality by means of  Littlewood-Paley decomposition and properties of singular homogeneous kernels of degree -$d$. A refinement when $q>2$ is proved.
\end{abstract}

The classical Hardy inequality states that when $d\geq 3$
\begin{equation}\label{eq:hardycl}
\int_{\R^d} \frac{|u|^2}{|x|^2}dx \leq \frac{4}{(d-2)^2}\int_{\R^d} |\nabla u|^2 dx
\end{equation}
and it is clearly of fundamental importance in analysis. There are of course many different  proofs of \eqref{eq:hardycl}, the simplest one consists in restrict by density to 
$D(\R^d\setminus \{0\})$, to observe that $\frac{1}{|x|^2}=-\frac 12 x\cdot \nabla (\frac{1}{|x|^2})$, then to  integrate by parts and eventually to apply Cauchy-Schwarz inequality.\\
A natural extension of \eqref{eq:hardycl} is in the framework of fractional Sobolev spaces $\dot H^s(\R^d)$. In this setting the following Hardy-type inequality holds
\begin{equation}\label{eq:hardyfr}
\int_{\R^d} \frac{|u|^2}{|x|^{2s}}dx\leq C ||f||_{\dot H^s (\R^d)}^2,
\end{equation}
provided that $0\leq s <\frac{d}{2}.$ For a compact and nice proof of \eqref{eq:hardyfr} we quote Theorem 2.57 in \cite{BCD}  and the proof given by Tao in the Appendix of \cite{T} while for an improvement involving Besov spaces we quote \cite{BCG}.\\
\\
If one is interested in proving an $L^q$ estimate for  $\frac{|f|}{|x|^s}$ we need to recall the definition of the  homogeneous Sobolev norm $||f||_{\dot W^{s,q}(\R^d)}$ which  is defined as $|||D|^s f||_{L^q(\R^d)}$ where 
$ (\widehat {|D|^s  f}) (\xi)= |2\pi\xi|^s \widehat{u} (\xi).$
In this note we give a direct proof and a refinement when $q>2$ for the following class of Hardy-type inequalities that generalize the fractional Hardy inequality \eqref{eq:hardyfr}.
\begin{thm}[Fractional Hardy inequality]\label{thm:hardysw}
Let  $0<s<\frac{d}{q}$, $1<q<\infty$ and $f\in \dot W^{s,q}(\R^d)$, then 
\begin{equation}\label{eq:hardy}
||\frac{f}{|x|^s}||_{L^q(\R^d)}\leq C(d,s,q)||f||_{\dot W^{s,q}(\R^d)}.
\end{equation}
\end{thm}
\noindent 
The explicit value of the constant $C(d,s,q)$ in \eqref{eq:hardy} is due to Herbst \cite{H}. The proof of \eqref{eq:hardy} goes back to the end of the fifties of the last century thanks to the work of Stein and Weiss \cite{SW} who proved 
an even more general version of  \eqref{eq:hardy} called Stein-Weiss inequality given by
\begin{equation}\label{eq:steinweiss}
\left(\int_{\R^d} \left(|T_{\lambda }f(x)| |x|^{-\beta} \right)^qdx\right)^{\frac{1}{q}} \leq C(d,q,p, \lambda) \left(\int_{\R^d }\left(|f(x)||x|^{\alpha}\right)^{p}dx\right)^{\frac{1}{p}}
\end{equation}
where 
$$T_{\lambda}f(x)=\int_{\R^d}\frac{f(y)}{|x-y|^\lambda}dy \ \ \ \ 0<\lambda<d,$$
and 
$$0<\lambda<d, 1<p<\infty, \alpha<\frac{d}{p'}, p\leq q<\infty, \  \beta<\frac{d}{q}, \ \alpha+\beta\geq 0,$$
$$\frac{1}{q}=\frac{1}{p}+(\frac{\lambda+\alpha+\beta}{d})-1.$$
The fact that \eqref{eq:steinweiss} implies \eqref{eq:hardy} follows by the fact that $T_{\lambda}f = c |D|^{-s} f, $ with $\lambda = d-s, $
$ c = \frac{\pi^{d / 2} \Gamma((d-\lambda) / 2)}{ \Gamma(\lambda / 2)}$
 and choosing $p=q$ and $\alpha=0, \beta=s$.\\
 \\
 
 In order to state our result we recall the standard definition for Homogeneous Besov norm $||\cdot ||_{\dot B_{p,q}^s}$ and Tribel-Lizorkin norm  $||\cdot ||_{\dot F_{p,q}^s}$ (see e.g. \cite{FJW} for general references).
 Let $f$ be a tempered distribution such that $\hat f\in L^1_{loc}$ and $P_N(f)$ the Littlewood-Paley projector on the dyadic frequency $N$, i.e. $\widehat{P_N(f)}(\xi)=\psi_N(\xi)\hat f(\xi)$ where $\psi_N(\xi)=\psi(\frac{\xi}{N})$ and $\sum_{N \in 2^{\Z}} \psi_N=1$, then we define
 $$||f||_{\dot B_{p,q}^s}=\left(\sum_{N\in 2^{\Z}} ||N^sP_N(f)||_{L^p}^q\right)^{\frac 1q},$$
 $$||f||_{\dot F_{p,q}^s}=||\left(\sum_{N\in 2^{\Z}} |N^sP_N(f)(x)|^q\right)^{\frac 1q}||_{L^p}.$$
 Our result is a direct proof of the following
 \begin{thm}\label{thm:hardyjm}
Let  $0<s<\frac{d}{q}$, $1<q<\infty$  then 
\begin{equation}\label{eq:hardyjm}
||\frac{f}{|x|^s}||_{L^q(\R^d)}\leq C(d,s,q)||f||_{\dot B^s_{q,q}(\R^d)},
\end{equation}
\end{thm}
with the following corollary
 \begin{cor}\label{cor:hardyjm}
Let  $0<s<\frac{d}{q}$, if $1<q\leq 2$ then   
\begin{equation}\label{eq:hardy2}
||\frac{f}{|x|^s}||_{L^q(\R^d)}\leq C(d,s,q)||f||_{\dot W^{s,q}(\R^d)},
\end{equation}
if  $q> 2$
\begin{equation}\label{eq:hardy3}
||\frac{f}{|x|^s}||_{L^q(\R^d)}\leq C(d,s,q)||f||_{\dot W^{s,q}(\R^d)}^{\frac{1}{q}}||f||_{\dot F^s_{q,2(q-1)}(\R^d)}^{\frac{q-1}{q}}.
\end{equation}
\end{cor}
The fact that $||\frac{f}{|x|^s}||_{L^q(\R^d)}$ can be controlled by homogeneous Besov norms is not a novely, a proof of Theorem \ref{thm:hardyjm} can be found in \cite{Y}, see also \cite{Y2}. Here we present a direct proof using the Shur test. We shall remark that  our corollary when $q>2$ is a refinement of Hardy inequality \eqref{eq:hardy}. Indeed we have  when $2(q-1)>2$
 $$||f||_{\dot F^s_{q,2(q-1)}(\R^d)}^{\frac{q-1}{q}}\leq ||f||_{\dot F^s_{q,2}(\R^d)}^{\frac{q-1}{q}}\sim ||f||_{\dot W^{s,q}(\R^d)}^{\frac{q-1}{q}}$$
 thanks to square function estimate
 $$||f||_{\dot F_{q,2}^s}=||\left(\sum_{N\in 2^{\Z}} |N^sP_N(f)(x)|^2\right)^{\frac 12}||_{L^q}\sim |||D|^s f||_{L^q(\R^d)}.$$
 The case $1<q<2$ is proved by duality and it requires proving the $L^q$ continuity for singular homogeneous kernels of degree -$d$. This fact is well known and is Lemma 2.1   in \cite{SW}. We underline however that our strategy in proving Theorem \ref{thm:hardyjm} permits to skip the more delicate lemmas in the Stein and Weiss paper \cite{SW} that are needed to prove \eqref{eq:hardy}.\\
 As a final comment, recalling that $|D|f=\sum_{j=1}^d R_j(\partial_{x_j} f)$ with $R_j$ the Riesz transform defined as $ (\widehat {R_j f}) (\xi)= -i \frac{\xi_j}{|\xi|} \widehat{u} (\xi)$
 and that hence $|||D|f||_{L^q(\R^d)}\lesssim ||\nabla f||_{L^q(\R^d)}$ when $1<q<\infty$, we get 
 \begin{cor}\label{cor:hardyjm2}
Let $2<q<d$ then   
\begin{equation}\label{eq:hardy4}
||\frac{f}{|x|}||_{L^q(\R^d)}\leq C(d,s,q)||\nabla f||_{L^q(\R^d)}^{\frac{1}{q}}||f||_{\dot F^s_{q,2(q-1)}(\R^d)}^{\frac{q-1}{q}}.
\end{equation}
\end{cor}
We underline that Corollary \ref{cor:hardyjm2} is a refinement of the classical Hardy inequality involving $\nabla f$
 \begin{equation}\label{eq:hardy5}
||\frac{f}{|x|}||_{L^q(\R^d)}\leq \left(\frac{q}{d-q}\right)||\nabla f||_{L^q(\R^d)}.
\end{equation}
 by the fact that $||f||_{\dot F^s_{q,2(q-1)}(\R^d)}\leq ||f||_{\dot F^s_{q,2}(\R^d)}\lesssim ||\nabla f||_{L^q(\R^d)}$.
In the literature there is a lot of interest in proving improvements for \eqref{eq:hardy5}, typically such improvement (in bounded or unbounded domains) are in the direction to add a negative term in r.h.s of \eqref{eq:hardy5}, see e.g. \cite{BV,BM,D, FT,FS,Ga, Gh, Ma}. Our refinement, although obtained with different techniques, is more in the spirit of  \cite{BCG} and \cite{Tr}, i.e. to control r.h.s. of \eqref{eq:hardy5} with terms that are smaller (up to a multiplicative constant) than the Sobolev norms.
 \section{Proof of Theorem \ref{thm:hardyjm}}
 A key argument in our proof is given by the following well known version of Shur test
\begin{prop}\label{shurtest}cLet $\alpha_{N,R}\geq 0$, with $N,R\in 2^{\Z}$, $1<q<\infty$, then
$$\sum_R\left(\sum_N \alpha_{N,R} C_N \right)^q\lesssim \sum_N \left(C_N\right)^{q}$$
provided there exists a sequence of positive numbers $p_N$ such that
\begin{eqnarray}
 \left(\sum_N \alpha_{N,R} p_N^{\frac{q'}{q}}\right)^{\frac{q}{q'}}\lesssim p_R\label{shur1}\\
 \sum_R \alpha_{N,R} p_R\lesssim p_N. \label{shur2}
\end{eqnarray}
\end{prop}
\begin{proof}

By Holder's inequality with conjugated exponent $(q, q')$
$$\sum_N \alpha_{N,R} C_N =\sum_N \alpha_{N,R}^{\frac{1}{q}}\alpha_{N,R}^{\frac{1}{q'}} p_N^{\frac{1}{q}}\frac{C_N}{p_N^{\frac{1}{q}}} \leq  \left(\sum_N \alpha_{N,R} p_N^{\frac{q'}{q}}\right)^{\frac{1}{q'}}  \left(\sum_N \alpha_{N,R} \frac{C_N^q}{p_N}\right)^{\frac{1}{q}}$$
we get
$$\sum_R\left(\sum_N \alpha_{N,R} C_N \right)^q\leq \sum_R \left(\sum_N \alpha_{N,R} p_N^{\frac{q'}{q}}\right)^{\frac{q}{q'}}  \left(\sum_N \alpha_{N,R} \frac{C_N^q}{p_N}\right)$$
that, thanks to \eqref{shur1} and Fubini, implies
$$\sum_R\left(\sum_N \alpha_{N,R} C_N \right)^q\lesssim \sum_R p_R  \left(\sum_N \alpha_{N,R} \frac{C_N^q}{p_N}\right)=\sum_N \frac{C_N^q}{p_N}\left(\sum_R \alpha_{N,R}p_R \right).$$
Now by \eqref{shur2} we conclude
$$\sum_R\left(\sum_N \alpha_{N,R} C_N \right)^q\lesssim \sum_N \frac{C_N^q}{p_N} p_N=\sum_N C_N^q.$$
\end{proof}

The strategy of the proof for  is an adaptation of  proof of Hardy inequality in the case $q=2$ given by Tao \cite{T}, i.e.  to prove the following estimate
\begin{equation}\label{assump1}
\int_{\R^d}\frac{|f(x)|^q}{|x|^{sq}}dx \lesssim  \sum_N N^{qs} ||P_N f||_{L^q(\R^d)}^q
\end{equation}
where $P_N f$ are the classical Littlewood-Paley projectors with $N$ a dyadic number. 

We devide $\R^d$ in dyadic shells obtaining

\begin{equation}
    \int_{\mathbb{R}^{d}}\frac{\vert f(x) \vert^{q}}{\vert x \vert^{qs}}\,dx = \sum_{R \in 2^{\mathbb{Z}}} \int_{\frac{R}{2}\leq \vert x \vert \leq R} \frac{\vert f(x)\vert^{q}}{\vert x\vert^{qs}}\,dx \lesssim \sum_{R \in 2^{\mathbb{Z}}} \frac{1}{R^{sq}}\,\int_{\{\frac{R}{2}\leq \vert x \vert\leq R\}} \vert f \vert^q \, dx. \label{primopasdisug:cap6}
\end{equation}
such that using the Littlewood-Paley decomposition we get
\begin{equation}
    \sum_{R \in 2^{\mathbb{Z}}} \frac{1}{R^{sq}}\,\int_{\{\frac{R}{2}\leq \vert x \vert\leq R\}} \vert f \vert^q \, dx\leq \sum_{R \in 2^{\mathbb{Z}}} R^{-sq}\, \left( \sum_{N \in 2^{\mathbb{Z}}}  \left( \int_{\{\vert x \vert \leq R\}} \vert P_{N}(f) \vert^{q} \right)^{\frac{1}{q}} \right)^{q}.
\end{equation}
By the Bernstein inequality $||P_{N}(f)||_{L^{\infty}(\R^d)}\leq N^{\frac{d}{q}}||P_{N}(f)||_{L^{q}(\R^d)}$
it follows that
\begin{equation}
    \left( \int_{\frac{R}{2}<\vert x \vert< R} \vert P_{N}(f) \vert^{q} \right)^{\frac{1}{q}} \leq R^{\frac{d}{q}}\Vert P_{N}(f) \Vert_{L^{\infty}} \\
    \leq (NR)^{\frac{d}{q}} \Vert P_{N}(f) \Vert_{L^{q}},
\end{equation}
and clearly 
\begin{equation*}
    \left( \int_{\frac{R}{2}<\vert x \vert< R} \vert P_{N}(f) \vert^{q} \right)^{\frac{1}{q}} \leq \Vert P_{N}f \Vert_{L^q},
   \end{equation*}
such that we get
$$
 \int_{\mathbb{R}^{d}}\dfrac{\vert f(x)\vert^{q}}{\vert x \vert^{qs}} dx\lesssim \sum_{R}R^{-qs} \left( \sum_{N} \min\{ 1,(NR)^{\frac{d}{q}}\}\, \Vert P_{N}f \Vert_{L^{q}} \right)^{q} \notag \\
=$$ 
$$=\sum_{R}\left(\sum_{N} \min\{ (NR)^{-s},(NR)^{\frac{d}{q}-s}\}\, \Vert \, N^{s}\,P_{N}f \Vert_{L^{q}}\right)^{q}.$$

The last step is to apply the Schur test given by Proposition \ref{shurtest} in order to conclude that
$$
    \sum_{R}\left(\sum_{N} \min\{ (NR)^{-s},(NR)^{\frac{d}{q}-s}\}\, \Vert \, N^{s}\,P_{N}f \Vert_{L^{q}}\right)^{q} \leq \sum_{N \in 2^{\mathbb{Z}}} N^{sq}\, \Vert P_{N}(f) \Vert_{L^{q}}^{q}\notag\\
    =$$
    $$=\sum_{N \in 2^{\mathbb{Z}}}N^{sq}\,\int_{\mathbb{R}^{d}}\vert P_{N}(f)\vert^{q}=\int_{\mathbb{R}^{d}} \sum_{N \in 2^{\mathbb{Z}}} N^{sq}\vert P_{N}(f)\vert^{q}. $$
Notice that 
\begin{eqnarray}
 \sum_{N>\frac{1}{R}}  \min\{ (NR)^{-s},(NR)^{\frac{d}{q}-s}\}+\ \sum_{N\leq \frac{1}{R}}  \min\{ (NR)^{-s},(NR)^{\frac{d}{q}-s}\}= \nonumber\\
= R^{-s}  \sum_{N>\frac{1}{R}}  N^{-s}+R^{\frac{d}{q}-s} \sum_{N\leq \frac{1}{R}} N^{\frac{d}{q}-s}\lesssim 1 \nonumber
\end{eqnarray}
such that (arguing in the same way when summing over $R$)
  \begin{eqnarray}
 \sum_N  \min\{ (NR)^{-s},(NR)^{\frac{d}{q}-s}\}  \lesssim 1\\
 \sum_R  \min\{ (NR)^{-s},(NR)^{\frac{d}{q}-s}\}  \lesssim 1. \end{eqnarray}
The hypoteses for Shur test  given by Proposition \ref{shurtest} are hence fulfilled by choosing $\alpha_{N,R}=  \min\{ (NR)^{-s},(NR)^{\frac{d}{q}-s}\} $ and $p_N=1$ in Proposition \ref{shurtest}.
 This proves \eqref{eq:hardy}. 
 \\

\section{Proof of Corollary \ref{cor:hardyjm}}
In Theorem \ref{thm:hardyjm} we proved the  following estimate
\begin{equation}\label{assump1}
\int_{\R^d}\frac{|f(x)|^q}{|x|^{sq}}dx \lesssim  \sum_N N^{qs} ||P_N f||_{L^q(\R^d)}^q
\end{equation}
where $P_N f$ are the classical Littlewood-Paley projectors with $N$ a dyadic number.  First we prove that \eqref{assump1} implies the Fractional Hardy inequality.
We have two cases: $q\geq 2, q<2$.\\
Case $q\geq 2$:\\

Thanks to  \eqref{assump1}
we derive
$$\sum_N N^{qs} ||P_N f||_{L^q(\R^d)}^q=\int_{\R^d} \sum_N N^{sq}|P_N f(x)|^q dx \leq  \int_{\R^d} \left(\sum |N^s P_N f(x)|^2\right)^{\frac{q}{2}} dx  $$
from the elementary inequality $\left(\sum_i a_i^{p_1}\right)^{\frac{1}{p_1}}\leq \left(\sum_i a_i^{p_2}\right)^{\frac{1}{p_2}}$ with $p_1\geq p_2,$
obtaining $$\int_{\R^d}\frac{|f(x)|^q}{|x|^{sq}}dx \lesssim \sum_N N^{qs} ||P_N f||_{L^q(\R^d)}^q\leq $$
$$ \leq  \int_{\R^d} \left(\sum_N |N^s P_N f(x)|^2\right)^{\frac{q}{2}} dx \sim |||D|^s f||_{L^q(\R^d)}^q$$
where the last equivalence is nothing but the classical square function estimate, see for instance \cite{MS}.\\
To prove \eqref{eq:hardy3} we notice that
$$\int_{\R^d} \sum_N N^{sq}|P_N f(x)|^q dx\leq $$
$$\leq \int_{\R^d} \left( \sum_N N^{2s}|P_N f(x)|^2 \right)^{\frac 12} \left( \sum_N N^{2s(q-1)}|P_N f(x)|^{2(q-1)} \right)^{\frac 12}dx \leq$$
$$\leq  \left(\int_{\R^d} \left( \sum_N N^{2s}|P_N f(x)|^2 \right)^{\frac q2}dx\right)^{\frac 1q}\left(\int_{\R^d} \left( \sum_N N^{2s(q-1)}|P_N f(x)|^{2(q-1)} \right)^{\frac{q}{2(q-1)}}dx\right)^{\frac{q-1}{q}}$$
by applying  twice the  Holder's inequality, first  in the serie with conjugated exponent $(2, 2)$ and then in the integral  with conjugated exponent $(q, \frac{q}{q-1})$.
By definition 
$$\left(\int_{\R^d} \left( \sum_N N^{2s(q-1)}|P_N f(x)|^{2(q-1)} \right)^{\frac{q}{2(q-1)}}dx\right)^{\frac{q-1}{q}}=||f||_{\dot F^s_{q, 2(q-1)}}^{q-1}.$$
\\
Case $q< 2$:\\
For the case $q <2$ we use the dual characterization of $L^q$ norms, i.e.

\begin{gather*}
    \Vert \frac{f}{\vert x \vert^{s}} \Vert_{L^{q}}=\sup_{\Vert g \Vert_{q'}=1} \langle \frac{ f(x)}{\vert x \vert^{s}},\,g \rangle= \sup_{\Vert g \Vert_{q'}=1}\, \langle f(x),\, \frac{g(x)}{\vert x \vert^{s}}\rangle\\ 
    =\sup_{\Vert g \Vert_{q'}=1}\, \langle \vert D \vert^{-s}(\vert D \vert^{s}  f(x)),\, \frac{g(x)}{\vert x \vert^{s}}\rangle= \sup_{\Vert g \Vert_{q'}=1} \langle \vert D \vert^{s}f ,\, \vert D \vert^{-s}(\frac{g(x)}{\vert x \vert^{s}})\rangle\\
    \leq \, \Vert \vert D \vert^{s}f \Vert_{L^{q}}\, \Vert \vert D \vert^{-s}(\frac{g(x)}{\vert x \vert^{s}})\Vert_{L^{q'}}.
\end{gather*}

Now we aim to prove that

\begin{equation}\label{eq:sw2}
   \Vert \vert D \vert^{-s}(\frac{g(x)}{\vert x \vert^{s}})\Vert_{L^{q'}(\R^d)} \lesssim \Vert g\Vert_{L^{q'}(\R^d)},
\end{equation}
for all  $g \in L^{q'}$ with $q'>2$ such that 
we could conclude that 
\begin{gather*}   
\Vert \frac{f}{\vert x \vert^{s}} \Vert_{L^{q}(\R^d)}=\sup_{\Vert g \Vert_{q'}=1} \langle \frac{\vert f(x)\vert}{\vert x \vert^{s}},\,g \rangle  \lesssim  \Vert D \vert^{s}f \Vert_{L^{q}(\R^d)}.
\end{gather*}

Now we prove \eqref{eq:sw2}. We have (skipping $q'$ with $q$ to simplify the notation)
\begin{gather*}
    \vert D \vert^{-s}(\frac{g(x)}{\vert x \vert^{s}})\vert^q  \sim \left\vert \int_{\mathbb{R}^{d}} \frac{g(y)}{\vert x-y\vert^{d-s}\,\vert y \vert^{s}}\,dy   \right\vert^{q} \leq \left\vert \int_{\mathbb{R}^{d}} \frac{\vert g(y)\vert}{\vert y \vert^{s}\, \vert x-y \vert^{d-s}} \,dy \right\vert^{q}\\
\lesssim \left\vert \int_{\mathbb{R}^{d}} \frac{\vert g(y)\vert\, \mathbbm{1}_{\{ \vert y \vert > \frac{\vert x\vert}{2}\}}(y) }{\vert y \vert^{s}\, \vert x-y\vert^{d-s}} dy\right\vert^{q}+\left\vert \int_{\mathbb{R}^{d}} \frac{\vert g(y)\vert\, \mathbbm{1}_{\{ \vert y \vert \leq \frac{\vert x\vert}{2}\}}(y) }{\vert y \vert^{s}\, \vert x-y\vert^{d-s}} dy\right\vert^{q}\\
\lesssim \frac{1}{\vert x \vert^{qs}}\, \left\vert \int_{\mathbb{R}^{d}} \frac{\vert g(y)\vert\, \mathbbm{1}_{\{ \vert y \vert > \frac{\vert x\vert}{2}\}}(y) }{ \vert x-y\vert^{d-s}} dy\right\vert^{q}+\left\vert \int_{\mathbb{R}^{d}} \frac{\vert g(y)\vert\, \mathbbm{1}_{\{ \vert y \vert \leq \frac{\vert x\vert}{2}\}}(y) }{\vert y \vert^{s}\, \vert x-y\vert^{d-s}} dy\right\vert^{q}\\
\lesssim \frac{1}{\vert x \vert^{qs}}\, \left\vert \int_{\mathbb{R}^{d}} \frac{\vert g(y)\vert}{ \vert x-y\vert^{d-s}}dy \right\vert^{q}+\left\vert \int_{\mathbb{R}^{d}} \frac{\vert g(y)\vert\, \mathbbm{1}_{\{ \vert y \vert \leq \frac{\vert x\vert}{2}\}}(y) }{\vert y \vert^{s}\, \vert x-y\vert^{d-s}}dy \right\vert^{q}\\
:=\vert S_{1}(g) \vert^{q}+\vert S_{2}(g)\vert^{q}
\end{gather*}

By  previous estimates using Paley-Littlewood decompostion and the square function equivalence we get 
when $q>2$
\begin{gather*}
    \int_{\mathbb{R}^{d}}\vert S_{1}(g)\vert^{q}\,dx\sim \int_{\mathbb{R}^{d}}\left\vert \frac{\vert D\vert^{-s}\vert g(x)\vert }{\vert x\vert^{s}} \right\vert^{q}\,dx \lesssim \left\Vert \vert D \vert^{s}(\vert D \vert^{-s}|g|) \right\Vert_{L^{q}(\R^d)}^q=\Vert g\Vert_{L^{q}(\R^d)}^q.
\end{gather*}

Concerning  $\Vert S_{2}(g)\Vert_{L^q}$ we follow the strategy of Stein and Weiss in \cite{SW}  proving the $L^q$ continuity for singular homogeneous kernels of degree -$d$. The proof of this fact is Lemma 2.1 in \cite{SW} that we show for reader convenience.
First notice that $\frac{\vert y \vert}{\vert x \vert}\leq \frac{1}{2}$ implies 
\begin{gather*}
    \vert x-y\vert \geq \vert x \vert - \vert  y \vert \geq \frac{\vert x \vert}{2},
\end{gather*}
such that
\begin{gather}
    \int_{\vert y \vert \leq \frac{\vert x \vert}{2}} \frac{\vert g(y)\vert}{\vert x-y\vert^{d-s}\vert y\vert^{s}}\, dy \lesssim \int_{\vert y \vert \leq \frac{\vert x \vert}{2}} \frac{\vert g(y)\vert}{\vert y \vert^{s}\, \vert x\vert^{d-s}}\,dy. \label{weiss:intdastim}
\end{gather}
Now we introduce following \cite{SW} the function, 
\begin{gather*}
    K(x,y)=\begin{cases}
    \vert y \vert^{s}\,\vert x \vert^{d-s} \quad \vert y \vert \leq \frac{\vert x \vert}{2}\\
    0 \quad \text{otherwise}
    \end{cases}
\end{gather*}
and 
\begin{gather*}
    Ug(x):=\int_{\vert y \vert \leq \frac{\vert x}{2} \vert} \frac{\vert g(y)\vert}{\vert y \vert^{s}\, \vert x\vert^{d-s}}\,dy=\int_{\mathbb{R}^{d}} K(\vert x\vert,\,\vert y \vert)|g(y)|\,dy.
\end{gather*}
To conclude the proof it suffices hence to show that
\begin{gather*}
    \int_{\mathbb{R}^{d}} \vert Ug \vert^{q}dx \lesssim \int \vert g \vert^{q} dx.
\end{gather*}
Fixing  $\eta \in S^{d-1}$ and calling  $\vert x\vert=R$ we define
\begin{gather*}
    U_{\eta}g(R):=\int_{0}^{+\infty} \,r^{d-1}\, K(R,r)\cdot |g(r\,\eta)|dr,
\end{gather*}
such that 
\begin{gather*}
Ug(x)=\int_{\mathbb{R}}K(\vert x\vert,\vert y \vert)|g(y)|dy=\int_{0}^{+\infty} \left(\int_{S^{d-1}} \hspace{-8pt} K(R,r)|g(r\,\eta)|\,d\sigma_{\eta}\right)\,r^{d-1}\,dr\\
=\int_{S^{d-1}}\int_{0}^{+\infty} \hspace{-8pt}K(R,r)|g(r\eta)|\,r^{d-1}\,dr\,d\sigma_{\eta}=\int_{S^{d-1}} U_{\eta}g(R)\,d\sigma_{\eta}.
\end{gather*}

By the substitution $r=tR$ we obtain 
\begin{gather*}
   \hspace{-3cm}U_{\eta}g(R) = \int_{0}^{+\infty}\hspace{-8pt} K(R,Rt)\,|g(t\,R\,\eta)|\,R^{d-1}\,t^{d-1}\,R\,dt\\
   =\int_{0}^{+\infty}\hspace{-8pt}K(1,t)\,|g(t\,R\,\eta)|t^{d-1}\,dt,
\end{gather*}
thanks to the fact that  $K$ is homogeneous of degree $-d$ , i.e. that 
\begin{gather*}
    K(\lambda x,\,\lambda y)=\vert \lambda \vert^{-d}K(\vert x\vert,\vert y\vert).
\end{gather*}
Let  $h$ be the function in  $L^{q'}((0,+\infty);R^{d-1}\,dR)$ of unitary norm  such that
\begin{gather*}
    \left( \int_{0}^{+\infty} \vert U_{\eta}g(R)\vert^{q}\,R^{d-1}\,dR \right)^{\frac{1}{q}}=\int_{0}^{+\infty} \hspace{-8pt}U_{\eta}g(R)h(R)\,R^{d-1}\,dR\\
    =\int_{0}^{+\infty} \left\{ \int_{0}^{+\infty}K(1,t)\,\vert g(t\,R\,\eta)\vert \,t^{d-1}\,dt \right\}\,R^{d-1}\,h(R)\,dR\\
    =\int_{0}^{+\infty} \hspace{-8pt} K(1,t)\,t^{d-1} \left\{ \int_{0}^{+\infty} \hspace{-8pt}\vert g(tR\eta) \vert\, h(R)\,R^{d-1}\,dR \right\}\,dt\\
    \leq \int_{0}^{+\infty}K(1,t)\,t^{d-1} \left\{\int_{0}^{+\infty}\vert g(t\,R\,\eta)\vert^{q}\,R^{d-1}\,dR  \right\}^{\frac{1}{q}}\,dt\\
    =\left( \int_{0}^{+\infty}\hspace{-8pt} K(1,t)t^{d-1-\frac{d}{q}}\,dt \right)\cdot \left\{ \int_{0}^{+\infty}\hspace{-8pt} \vert g(R\eta)\vert^{q}\,R^{d-1}\,dR\right\}^{\frac{1}{q}}\\
    =\left(\int_{0}^{1}t^{d-\frac{d}{q}-1-s}\,dt\right)\cdot \left\{ \int_{0}^{+\infty} \hspace{-8pt} \vert g(R\eta)\vert^{q}\,R^{d-1}\,dR\right\}^{\frac{1}{q}}=:J\cdot \left\{ \int_{0}^{+\infty} \hspace{-8pt} \vert g(R\eta)\vert^{q}\,R^{d-1}\,dR\right\}^{\frac{1}{q}} ,
\end{gather*}
where the last integral  $J$ converges due to the fact that by our assumptions  $s <\frac{d}{q'}$ (remember that we skipped $q'$ with $q$).

Now we estimate $L^{q}(\mathbb{R}^{d})$ norm of $Ug. $ By Jensen inequality
\begin{gather*}
    \vert Ug(R) \vert^{q}=\left\vert \int_{S^{d-1}} \vert U_{\eta}g(R)\vert\,d\sigma_{\eta} \right\vert^{q} \leq \{\vert S^{d-1}\vert \}^{q-1}\,\int_{S^{d-1}} \vert U_{\eta}g\vert^{q}\,d\sigma_{\eta},
\end{gather*}
such that integrating with respect to the measure  $R^{d-1}dR$ we get 
\begin{gather*}
    \hspace{-4cm}\int_{0}^{+\infty} \hspace{-5pt}\vert Ug(R) \vert^{q}R^{d-1}\,dR
    \leq \\
    \leq J^{q}\,\vert S^{d-1} \vert^{q-1}\left( \int_{0}^{+\infty} \left\{ \int_{S^{d-1}} \hspace{-5pt} \vert U_{\eta}g(R)\vert^{q}\,d\sigma_{\eta} \right\}\,R^{d-1}\,dR \right)\\
    =J^{q}\,\vert S^{d-1} \vert^{q-1}\,\int_{S^{d-1}}\int_{0}^{+\infty} \vert U_{\eta}g(R) \vert^{q}\,R^{d-1}\,dR \, d\sigma_{\eta}\\
    \leq J^{q}\,\vert S^{d-1} \vert^{q-1}\, \int_{S^{d-1}}\int_{0}^{+\infty}\vert g(R\,\eta) \vert^{q}\,R^{d-1}\,dR\,d\sigma=J^{q}\,\vert S^{d-1}\vert^{q-1}\,\int_{\mathbb{R}^{d}}\vert g(x)\vert^{q}\,dx.  
\end{gather*}
By the fact that  $Uf(x)$ is radial we can conclude that
\begin{gather*}
    \int_{\mathbb{R}^{d}} \vert Ug(x) \vert^{q}\,dx =\vert S^{d-1} \vert\cdot \int_{0}^{+\infty}\vert Ug(R)\vert^{q}\,R^{d-1}\,dR \leq J^{q}\vert S^{d-1}\vert^{q}\,\int \vert g(x)\vert^{q}\,dx.
\end{gather*}
This concludes the proof in the case $q<2$.
\\
\\
\textbf{DECLARATIONS}\\
\\
\noindent \textbf{Conflict of interest}. The authors declare that they have no conflict of interest.

\end{document}